\documentclass[11pt]{article}
\usepackage{amssymb}   
\usepackage{amsthm}    
\usepackage{amsmath}   
\usepackage{stmaryrd}  
\usepackage{titletoc}  
\usepackage{mathrsfs}  
\usepackage{cases}
\usepackage{color}
\newenvironment{enumerate-roman}{\begin{enumerate}}{\end{enumerate}}

\theoremstyle{plain}

\numberwithin{equation}{section} 
\title{\huge A Stochastic Generalized Ginzburg-Landau Equation Driven by Jump Noise
}
\author{ \Large Lin Lin\footnote{Jiangsu Provincial Key Laboratory for NSLSCS, School of Mathematical Sciences,  Nanjing Normal University, Nanjing 210023, China. (Email: \texttt{linlin@njnu.edu.cn})
}
~~~and~~~\Large Hongjun Gao\footnote{Corresponding author.  School of Mathematical Sciences and Jiangsu Center for Collaborative Innovation in Geographical Information Resource Development and Application,  Nanjing Normal University, Nanjing 210023, China. (Email:
\texttt{gaohj@njnu.edu.cn, hjgao-nj@163.com})}
}
\date{}
\allowdisplaybreaks
\begin{document}
\maketitle

 \textbf{Abstract} This paper is concerned with the stochastic generalized Ginzburg-Landau equation driven by a multiplicative
noise of jump type. By a prior estimate, weak convergence and monotonicity technique, we prove the existence and uniqueness of the solution of an initial-boundary value problem with homogeneous Dirichlet boundary condition. However, for the generalized Ginzburg-Landau equation, such a locally monotonic condition of the nonlinear term can not be satisfied in a straight way. For this, we utilize the characteristic structure of nonlinear term and refined analysis to overcome this gap.



\large
\section{Introduction}

The deterministic Ginzburg-Landau partial differential
equation had been used to model phenomena in a number of
 different areas in physics and other fields \cite{GH,GV,TR}, and it was extended to the generalized Ginzburg-Landau equation with derivative nonlinear term by Doelman in \cite{DA}. Many results on the existence and uniqueness were studied under various assumptions on the parameters \cite{GB,GD1,GD2,GH,GV,LG}. However, some perturbations may be neglected in the derivation of this ideal model. Researchers often represent the microscopic effects by random perturbations in the dynamics of the macroscopic observables. Thus, it is natural to consider stochastic effect in the Ginzburg-Landau equation. Recently, the stochastic Ginzburg-Landau equation with additive or multiplicative Gaussian noise has been studied by a few authors (see, e.g.\cite{WG,Y1,Y2,Y3}), among which Yang \cite{Y2,Y3} considered large deviations for the generalized $1$-D stochastic Ginzburg-Landau equation and existence for the generalized $2$-D stochastic Ginzburg-Landau equation with multiplicative Gaussian noise.

  As it is well known, most of the works of the stochastic generalized Ginzburg-Landau equation are driven by Gaussian white noise, but stochastic partial differential equations driven by jump noise have important applications in mathematical physics \cite{AR,SZ}. The applied backgrounds of jump noise is also our main motivation to consider the Ginzburg-Landau equation driven by jump noise. There are lots of existing results related to stochastic partial differential equations driven by jump noise, such as \cite{BZ,BH,BL,DX1,DX2,LR,RZ,SG}. To be specific, Brze\'{z}niak and Zhu \cite{BZ} obtained the existence and uniqueness of solutions for a type of stochastic nonlinear beam equations with Poisson-type noises; Brze\'{z}niak \cite{BH} considered the $2$-D stochastic Navier-Stokes equations driven by jump noise; and Sun and Gao \cite{SG} studied the well-posedness for $2$-D stochastic Primitive equations with L\'{e}vy noise.

In Brze\'{z}niak and Liu \cite{BL}, the authors built a unified framework for SPDE with locally monotonic coefficients driven by L\'{e}vy noise, in which the existence and uniqueness of the solution in a fixed probability space is proved based on a prior estimate, weak convergence and monotonicity arguments as in \cite{PR}. However, this method does not work for the generalized Ginzburg-Landau equation since the local monotonicity of the nonlinear operator cannot be satisfied in a straight way.  In order to obtain the local monotonicity in (\ref{27}), we first use the nonlinear structure and refined analysis (such as the analysis of $\int_0^T\int_D|u_n|^{2\sigma}|\nabla u_n|^2dxdt$) to overcome the gaps. Also the uniform bounds for $u_n$ in $L^p(\Omega, L^{\infty}([0,T], V))$ is established for some $p \geq 2$ (Lemma 3.3), which is important to enable $r(t)$ in (\ref{26}) to make sense. Moreover, different from the deterministic case in which $\textmd{Re}\, I\leq 0$ is enough, the more elaborate estimate $\textmd{Re} \,I\leq -C\|u_n-\phi\|_{2\sigma+2}^{2\sigma+2}$, $C>0$ is needed (Lemma 3.5), which is helpful for the estimates of some items of $\textmd{Re} \,J$ in Lemma 3.6 ($I,\ J$ are defined in Section 3).\\
\indent In this paper, we will study the following stochastic generalized Ginzburg-Landau equation (SGGLE) driven by jump noise:
\begin{numcases}{}\label{1}
du=((1+i\alpha)\triangle u-(1-i\beta)|u|^{2\sigma}u +\gamma u+F(u)) dt+ \int_Z g(t,u,z)\tilde{\eta}(dz,dt), \nonumber\\
u(x,t)=0,\hspace{6cm} t\geq 0,\,x\in\partial D,\\
u(x,0)=u_0(x),\hspace{6.6cm}  x\in D,\nonumber
\end{numcases}
where $D=(0,L_1)\times(0,L_2)$, $i=\sqrt{-1}$, $\gamma>0$, the
parameters $\alpha, \,\beta$ are all real-valued constants
and  $u$ is a complex-valued scalar
function. The derivative term $F(u)=\lambda_1\cdot\triangledown(|u|^2u)+(\lambda_2\cdot\triangledown u)|u|^2$ with two complex constant vectors $\lambda_1$ and $\lambda_2$. $\tilde{\eta}$ is the L\'{e}vy process defined on a complete probability space and $g$ is a given function which will be defined later. For the deterministic $2$-D generalized Ginzburg-Landau equation, Li and Guo \cite{LG} proved that the equation has a unique solution under the following assumptions on the parameters $\alpha, \beta$ and $\sigma$:

(1) Either (a) $\sigma>2$ or (b) $\sigma=2$ and $|\lambda_i|$, $i=1,2$, are sufficiently small;

(2) $-1+\alpha\beta<\frac{\sqrt{2\sigma+1}}{\sigma}|\alpha+\beta|$.

The next two sections are organized as follows. In Section 2, we recall some fundamental concepts related to L\'{e}vy process and present the main theorem. In section 3, we use the Galerkin method to prove the main theorem, i.e. existence and uniqueness of solution to our concerned equation.\\

\section{ Preliminaries and main theorem}

  Firstly, we introduce some definitions and basic properties of L\'{e}vy processes. The readers can also refer to \cite{PZ} for more details.

   Let $(\Omega ,\mathscr{F},\mathbb{F},\mathbb{P})$ be a filtered probability space, where $\mathbb{F}=(\mathscr{F}_t)_{t\geq 0}$ is a filtration, $(Z,\mathcal{Z})$ be a measurable space, and $\nu$ be a $\sigma$-finite positive measure on it. We denote the Borel $\sigma$-field on a topological space X by $\mathcal{B}(X)$. Let $\eta: \mathscr{F}\times\mathcal{B}(R^{+})\times\mathcal{Z}\rightarrow \bar{\mathbb{N}}=\mathbb{N}\cup\{\infty\}$ be a time homogeneous Poisson random measure with the intensity measure $\nu$ defined over $(\Omega ,\mathscr{F},\mathbb{F},\mathbb{P})$. We denote by $\tilde{\eta}\,(dt,dz)=\eta\,(dt,dz)-dt\,\nu(dz)$ the compensated Poisson random measure associated to $\eta$.

Suppose that $(H, |\cdot|_{H})$ is a Hilbert space. We then define a unique continuous linear operator $I$ serving as the stochastic integration with respect to the $\mathbb{F}$-predictable process $\xi: [0,T]\times Z\times \Omega \rightarrow H$ with \begin{eqnarray}\label{2}\mathbb{E}\int_0^T\int_Z|\xi(r,z)|_{H}^2\nu(dz)dr<\infty,\,\,\,\,T>0.\end{eqnarray}
For any random step process $\xi$ satisfying the condition (\ref{2}) with a representation$$\xi(r)=\sum\limits_{j=1}^{n}1_{(t_{j-1},t_{j}]}(r)\xi_{j},\,\,\,r\geq 0,$$where $\{0=t_0<t_1<...<t_n<\infty \}$ is a partition of $[0,\infty)$ and $\xi_j$ is an $\mathscr{F}_{t_{j-1}}$ measurable random variable for $j=1,2,\cdots, n$, we define $I(\xi)$ to be an $H$-valued adapted and c\`{a}dl\`{a}g process as follows \begin{eqnarray}I(\xi)(t)=\sum\limits_{j=1}^{n}\int_Z\xi_j(z)\tilde{\eta}(dz,(t_{j-1}\wedge t,t_j\wedge t]),\ \ \ t\geq 0.\nonumber\end{eqnarray}
Usually, we write \begin{eqnarray}I(\xi)(t):=\int_0^t\int_Z\xi(r,z)\tilde{\eta}(dr,dz),\,\,\,\,t\geq 0.\nonumber\end{eqnarray}
The continuity of the operator $I$ means that \begin{eqnarray}\mathbb{E}|\int_0^t\int_Z\xi(r,z)\tilde{\eta}(dr,dz)|_H^2=\mathbb{E}\int_0^t\int_Z|\xi(r,z)|_H^2\nu(dz)dr,\,\,\,\,t\geq 0. \nonumber\end{eqnarray}
For fixed $T>0$, we denote by $\mathcal{M}\,^2(0,T;L^2(Z,\nu;H))$ the class of all $\mathbb{F}$-predictable processes $\xi:\,[0,T]\times Z\times \Omega \rightarrow H $ satisfying the condition (\ref{2}), where $L^2(Z,\nu;H)$ is defined as the class of all functions $\eta:\,Z\rightarrow H$ satisfying $\int_Z|\eta(z)|_{H}^2\nu(dz)<\infty$.

Before giving the proof of the existence and uniqueness of solution, we need clarify some useful lemmas.\\

It\^{o} formula in Hilbert space is needed in our situation and one can refer to It\^{o} formula for semimartingale in \cite{me} (Theorem 27.2) for $\phi: E \rightarrow G$ be a $C^2$ function or the following one for the process given by (\ref{linlin}) below in \cite{zh} (Theorem 3.5.3).\\
{\bf Lemma 2.1} Assume that $E$ is an Hilbert space. Let $X$ be a process given by
 \begin{eqnarray}\label{linlin}X_t=X_0+\int_0^t a(s)ds+\int_0^t\int_Z f(s,z)\tilde{\eta}(ds,dz),\,\,\ \ t\geq 0, \end{eqnarray}
 where $a$ is an $E$-valued progressively measurable process on the space $\big(\mathbb{R}_+\times\Omega,\ \mathcal{B}_{\mathbb{R}_+}\times\mathscr{F}\big)$ such that for all $t\geq 0$, $\int_0^t \|a(s,w)\|_E ds<\infty$ a.s. and $f$ is a predictable process on $E$ with $\mathbb{E}\int_0^T\int_Z|f(s,z)|^2\nu(dz)ds<\infty$. Let $G$ be a separable Hilbert space and $\phi: E \rightarrow G$ be a function of $C\,^1$ such that $\phi\,^{'}$ is H\"{o}lder continuous. Then for each $t>0$, we have $\mathbb{P}$-a.s.
\begin{eqnarray*}\phi(X_t)&=&\phi(X_0)+\int_0^t\phi\,^{'}(X_s)a(s)ds+\int_0^t\int_Z[\phi\,^{'}(X_{s-})f(s,z)]\tilde{\eta}(ds,dz)\\
&&+\int_0^t\int_Z[\phi(X_{s-}+f(s,z))-\phi(X_{s-})-\phi\,^{'}(X_{s-})f(s,z)]\eta(dzds).\end{eqnarray*}\\
{\bf Lemma 2.2 (Burkholder-Davis-Gundy Inequality \cite{I})} For any $p\geq 2$ and $T>0$, there exists a constant $C_p$ such that for any real-valued square integrable c\`{a}dl\`{a}g martingale $M$ with $M(0)=0$, we have
 \begin{eqnarray}
\sup\limits_{0\leq t\leq T}|M(t)|^p\leq C_p \mathbb{E}[M]_T^{p/2},\nonumber\end{eqnarray}
where $[M]_t$, $0\leq t\leq T$, is the Meyer process of $M$.\\
{\bf Lemma 2.3 (Okazawa and Yokota \cite{OY})}
Let $H$ be a complex Hilbert space with inner product
$\langle\cdot,\cdot \rangle$ and norm $\|\cdot\|_H$. For $p\in(1,\infty)$
and any non-zero $z,w\in H$ with $z\neq w$, we have the
following inequality:
\begin{eqnarray}\frac{|\,\textmd{Im}\langle\|z\|_H^{p-2}z-\|w\|_H^{p-2}w,z-w\rangle|}{\textmd{Re}\langle\|z\|_H^{p-2}z-\|w\|_H^{p-2}w,z-w\rangle}
\leq\frac{|p-2|}{2\sqrt{p-1}}.\nonumber\end{eqnarray}

\indent For the mathematical setting of our problem, we introduce complex Sobolev spaces. Denote $(\cdot,\cdot)$ the inner product and the norm in $L^2(D)$, where
$$(u,v)=\textmd{Re}\int_D u(x)\bar{v}(x)dx,$$ for $u,v\in L^2(D)$.\\
\indent We always write $H=L^2,\,V=H_0^1,\,\|\cdot\|_p=\|\cdot\|_{L^p}$ and $\|\cdot\|=\|\cdot\|_2$ when $p=2$ for simplicity. Let $Au=(1+i\alpha)\triangle u $ and $Bu=-(1-i\beta)|u|^{2\sigma}u  +\gamma u+F(u)$. The operator $A$ is an isomorphism from $D(A)=V\cap H^2$ onto $H$. We now write (\ref{1}) in the following abstract form:
\begin{numcases}{}\label{3}
du(t)=(Au(t)+Bu(t))dt+ \int_Z g(t,u,z)\tilde{\eta}(dz,dt), \nonumber\\
u(x,t)=0,\hspace{3cm} t\geq 0,\,x\in\partial D,\\
u(x,0)=u_0(x),\hspace{2.5cm}  x\in D.\nonumber
\end{numcases}
To obtain the existence of solution to (\ref{3}), we assume that the function $g:[0,\infty)\times H \rightarrow L^2(Z,\nu;H)$ is a measurable function and  there exist nonnegative constants $k_i,i=1,2,3,4$ , $h\in L^p(Z,\nu;H)$ $(2\leq p<2\sigma)$ such that for any $t\in [\,0,T]$ and all $u,v\in V,$  \\
\indent $(C_1): \,\|g(t,u)\|_{L^2(Z,\nu;H)}^2\leq k_1\|u\|^2+k_2\|\nabla u\|^2;$\\
\indent $(C_2): \,\|g(t,u)-g(t,v)\|_{L^2(Z,\nu;H)}^2\leq k_3\|u-v\|^2+k_4\|\nabla(u-v)\|^2;$\\
\indent $(C_3):\,\,|\partial_u g(t,u,z)|\leq h(z)|u|$ for all $z\in Z.$\\
{\bf Definition 2.4 }An $H$-valued c\`{a}dl\`{a}g $\mathbb{F}$-adapted process $u(t)$ is said to be a solution of the stochastic generalized Ginzburg-Landau equation with jump noise (\ref{3}) if for its $dt\times d\mathbb{P}$-equivalent class $\tilde{u}$ we have\\
 \indent(1) $\tilde{u}\in L^p(\Omega;L^{\infty}([0,T];\,V))$; \\
\indent (2) For every $t\in[0,T]$, the following equality holds  $\mathbb{P}$-a.s.:\\
\begin{eqnarray}u(t)=u_0+\int_0^tA\tilde{u}(s)ds+\int_0^tB\tilde{u}(s)ds+\int_0^t\int_Z g(s,\tilde{u}(s),z)\tilde{\eta}(dz,ds).\nonumber\end{eqnarray}
\indent The main result of this paper is  the following theorem.\\
{\bf Theorem 2.5} Suppose that conditions $(C_1)-(C_3)$ hold, $ 2\leq p<2\sigma$, $\sigma>2$ and $0<|\beta|< \frac{\sqrt{2\sigma+1}}{\sigma}.$ Then for any $H$-valued function $u_0$ satisfying $\mathbb{E} \|\nabla u_0\|^{p}<\infty$, SPDE (\ref{3}) has a unique solution $u=u(t)$, $0\leq t \leq T$, satisfying
\begin{eqnarray}\mathbb{E}&\big(\sup\limits_{0\leq s\leq T}\|u(s)\|^2+\int_0^T\|\nabla u(s)\|^{2}ds+\int_0^T\|u(s)\|_{2\sigma+2}^{2\sigma+2}ds\big)\leq C\,(\mathbb{E}\|u_0\|^2+1)\nonumber\end{eqnarray}
and the additional regularity $$u(t)\in L^2(\Omega;L^2([0,T]; H^2))\cap L^p(\Omega;L^\infty([0,T];\,V)).$$
\section{Proof of main theorem}
The main method for the proof of Theorem 2.5 is the Galerkin approximation of (\ref{3}) and we divide the proof into four steps. The first three steps give the existence proof of solution to (\ref{3}) and the proof of uniqueness of solution is demonstrated in step 4.

 Step 1: Suppose that $\{e_i:i\in N\}\subset D(A)$ is an orthonormal basis of $H$ such that ${\rm span}\{e_i:i\in N\}$ is dense in $V$. Denote $H_n:={\rm span}\{e_1,...e_n\}$. Set $P_n:H\rightarrow H_n$ to be \\
$$P_nx=\sum\limits_{i=1}^{n}\langle x,e_i\rangle \,e_i,$$ where $P_n$ is the orthogonal projection onto $H_n$ in $H$.

 For simplicity, we denote $G(u)=A(u)+B(u)$. For each finite $n\in N$, we consider the following equation:
\begin{eqnarray}\label{4}
u_n(t)=u_n(0)+\int_0^tP_nG(u_n(s))ds+\int_0^t\int_Z P_ng(s,u_n(s),z)\tilde{\eta}(dz,ds),\ \  t\in[0,T],
\end{eqnarray}
where $u_n(0)=P_nu_0$. According to \cite{AB}, (\ref{4}) has a unique c\`{a}dl\`{a}g strong solution.

Then we give some priori estimates as preparation of proof of Theorem 2.5.\\
{\bf Lemma 3.1} Under the same assumptions as in Theorem 2.5, there exists a constant $C$ such that
\begin{eqnarray}\sup\limits_{n}\mathbb{E}\left(\sup\limits_{0\leq t\leq T}\|u_n(t)\|^2+\int_0^T\|\nabla u_n(s)\|^2ds+\int_0^T\|u_n(s)\|_{2\sigma+2}^{2\sigma+2}ds\right)\leq C(\mathbb{E}\|u_0\|^2+1).\nonumber\end{eqnarray}
$Proof:$ Applying It\^{o} formula to the process $\|u_n(t)\|^2$ and taking the real part, we obtain
\begin{eqnarray}\label{7}&&\|u_n(t)\|^2\nonumber\\
&=&\|u_n(0)\|^2+2\,\textmd{Re}\int_0^t (u_n(s),\,(1+i\alpha)\triangle u_n(s))ds+2\int_0^t(u_n(s),\gamma u_n(s))ds\nonumber\\&&
-2\,\textmd{Re}\int_0^t(u_n(s),(1-i\beta)|u_n(s)|^{2\sigma}u_n(s)+F(u_n(s)))ds\nonumber\\&&+\int_0^t\int_Z \|P_ng(s,u_n(s),z)\|^2\eta(ds,dz)+2\,\textmd{Re}
\int_0^t\int_Z(u_n(s-),P_ng(s,u_n(s),z))\tilde{\eta}(ds,dz)\nonumber\\&=&\|u_n(0)\|^2-2\,\int_0^t\|\nabla u_n(s)\|^2ds+2\,\gamma\int_0^t\|u_n(s)\|^2ds\nonumber\\&&-2\,\int_0^t\|u_n(s)\|_{2\sigma+2}^{2\sigma+2}ds+\,I_1(t)+\,I_2(t)+\,I_3(t),\end{eqnarray}
where \begin{eqnarray*} I_1(t)&=&\int_0^t\int_Z \|P_ng(s,u_n(s),z)\|^2\eta(ds,dz),\\
 I_2(t)&=& 2\,\textmd{Re}\int_0^t\int_Z(u_n(s-),P_ng(s,u_n(s),z))\tilde{\eta}(ds,dz),\\
 I_3(t)&=& -2\,\textmd{Re}\int_0^t(u_n(s),\lambda_1\cdot\triangledown(|u_n(s)|^2u_n(s))+(\lambda_2\cdot\triangledown u_n(s))|u_n(s)|^2)ds.\end{eqnarray*}
Now, for each natural number $R$, we consider the stopping time $\tau_R^n:=\inf\{t\geq 0:\|u_n(t)\|^2 \geq R\}\wedge T$. Since the process $u_n(t), t\in[0,T]$, is adapted and  c\`{a}dl\`{a}g, it is obvious that $\tau_R^n\uparrow T$, and $\mathbb{P}\{\tau_R^n<T\}=0$ as $R\rightarrow \infty$.

First we have
\begin{eqnarray}\label{6}\mathbb{E}\sup\limits_{0\leq s\leq t\wedge \tau_R^n}I_1(s)&=&\mathbb{E}\int_0^{t\wedge \tau_R^n}\int_Z \|P_ng(s,u_n(s),z)\|^2\nu(dz)ds\nonumber\\&\leq&k_1\,\mathbb{E}\int_0^{t\wedge \tau_R^n}\|u_n(s)\|^2ds+k_2\,\mathbb{E}\int_0^{t\wedge \tau_R^n}\|\nabla u_n(s)\|^2ds\end{eqnarray}
Here the first equality is due to the martingale property of the stochastic integration with repect to the compensated Poisson random measure $\tilde{\eta}\,(dt,dz)=\eta\,(dt,dz)-dt\,\nu(dz)$ and the second inequality is due to the application of ($C_1$).

And the process $I_2$ is a martingale, we can apply the B-D-G inequality in Lemma 2.2 and the condition ($C_1$) again to get
\begin{eqnarray}\mathbb{E}\sup\limits_{0\leq s\leq t\wedge \tau_R^n}I_2(s)&\leq& 2\,\mathbb{E} \sup\limits_{0\leq s\leq t\wedge \tau_R^n}\left|\int_0^s\int_Z(u_n(r-),P_ng(s,u_n(s),z))\tilde{\eta}(dr,dz)\right|\nonumber\\&\leq&
C_1\,{\mathbb{E}}\left[\int_0^{t\wedge \tau_R^n}\int_Z\|u_n(s)\|^2\|P_ng(s,u_n(s),z)\|^2\nu(dz)ds\right]^\frac{1}{2}\nonumber\\
&\leq& C_1\left[2C_1\,\mathbb{E}\int_0^{t\wedge \tau_R^n}(k_1\|u_n(s)\|^2+k_2\|\nabla u_n(s)\|^2)ds\right]^{\frac{1}{2}}\times\nonumber\\&&\left[\frac{1}{2C_1}\,\mathbb{E}\sup\limits_{0\leq s\leq t\wedge \tau_R^n}\|u_n(s)\|^2\right]^{\frac{1}{2}}\nonumber\\&\leq&\frac{1}{2}\,\mathbb{E}\sup\limits_{0\leq s\leq t\wedge \tau_R^n}\|u_n(s)\|^2+2C_1\,^{2}k_1\mathbb{E}\int_0^{t\wedge \tau_R^n}\|u_n(s)\|^2ds\nonumber\\&&+2C_1\,^2k_2\mathbb{E}\int_0^{t\wedge \tau_R^n}\|\nabla u_n(s)\|^2ds.\end{eqnarray}

Then we estimate the term $I_3$.  Since
\begin{eqnarray}\label{9}\lambda_1\cdot\nabla(|u_n|^2u_n)+(\lambda_2\cdot\nabla u_n)|u_n|^2\emph{}&=&((2\lambda_1+\lambda_2)\cdot\nabla u_n)|u_n|^2+(\lambda_1\cdot\nabla u_n)u_n^2,\end{eqnarray} we have
\begin{eqnarray}\label{5}\mathbb{E}\sup\limits_{0\leq s\leq t\wedge \tau_R^n}I_3(s)&\leq& (6|\lambda_1|+2|\lambda_2|)\,\mathbb{E} \sup\limits_{0\leq s\leq t\wedge \tau_R^n}
\int_0^s\int_D|u_n(s,x)|^3|\nabla u_n(s,x)|dxds\nonumber\\&\leq&(6|\lambda_1|+2|\lambda_2|)\,\mathbb{E}\int_0^{t\wedge \tau_R^n}\|u_n(s)\|_6^3\|\nabla u_n(s)\|ds
\nonumber\\&\leq&\frac{1}{2}\,\mathbb{E}\int_0^{t\wedge \tau_R^n}\|\nabla u_n(s)\|^2ds+2(3|\lambda_1|+|\lambda_2|)^2\,\mathbb{E}\int_0^{t\wedge \tau_R^n}\|u_n(s)\|_6^6ds
\nonumber\\&\leq&\frac{1}{2}\,\mathbb{E}\int_0^{t\wedge \tau_R^n}\|\nabla u_n(s)\|^2ds+\frac{1}{2}\,\mathbb{E}\int_0^{t\wedge \tau_R^n}\| u_n(s)\|_{2\sigma+2}^{2\sigma+2}ds\nonumber\\
&&+ C_2(|\lambda_1|,|\lambda_2|)\mathbb{E}\int_0^{t\wedge \tau_R^n}\|u_n(s)\|^2ds.
\end{eqnarray}
In (\ref{5}), we use the Young inequality in the last inequality as below
\begin{eqnarray}\label{20}\|u_n\|_6^6\leq\|u_n\|_{2\sigma+2}^\theta\|u_n\|^{1-\theta}\leq \epsilon\|u_n\|_{2\sigma+2}^{2\sigma+2}+c(\epsilon)\|u_n\|^2.\end{eqnarray}
Putting (\ref{6})-(\ref{5}) into (\ref{7}), we have
\begin{eqnarray}&&\frac{1}{2}\,\mathbb{E}\sup\limits_{0\leq s\leq t\wedge \tau_R^n}\|u_n(s)\|^2+\frac{C_3}{2}\,\mathbb{E}\int_0^{t\wedge \tau_R^n}\|\nabla u_n(s)\|^2ds+\frac{3}{2} \mathbb{E}\int_0^{t\wedge \tau_R^n}\|u_n(s)\|_{2\sigma+2}^{2\sigma+2}ds\nonumber\\&\leq &\mathbb{E}\|u_n(0)\|^2+\frac{C_4}{2}\,\mathbb{E}\int_0^{t\wedge \tau_R^n}\|u_n(s)\|^2ds\nonumber,\end{eqnarray}
where $\frac{C_3}{2}=\frac{3}{2}-(2C_1^2+1)k_2$, $\frac{C_4}{2}=2\gamma +(2C_1^2+1)k_1+C_2(|\lambda_1|,|\lambda_2|)$. Therefore, if $k_2$ is small enough such that $\frac{3}{2}-(2C_1^2+1)k_2>0$, we may apply Gronwall lemma to have
\begin{eqnarray}
&&\mathbb{E}\sup\limits_{0\leq s\leq t\wedge \tau_R^n}\|u_n(s)\|^2+\mathbb{E}\int_0^{t\wedge \tau_R^n}\|\nabla u_n(s)\|^2ds+\mathbb{E}\int_0^{t\wedge \tau_R^n}\|u_n(s)\|_{2\sigma+2}^{2\sigma+2}ds\nonumber\\&\leq& C(\mathbb{E}\|u_n(0)\|^2+1)\nonumber.
\end{eqnarray}
Recall that $\tau_R^n\uparrow T$, and $\mathbb{P}\{\tau_R^n<T\}=0$ as $R\rightarrow \infty$. It follows from Fatou lemma that
\begin{eqnarray*}&&\mathbb{E}\sup\limits_{0\leq t\leq T}\|u_n(t)\|^2+\mathbb{E}\int_0^T\|\nabla u_n(s)\|^2ds+\mathbb{E}\int_0^T\|u_n(s)\|_{2\sigma+2}^{2\sigma+2}ds\\&\leq&\liminf\limits_{n\rightarrow\infty}\mathbb{E}\left(\sup\limits_{0\leq s\leq t\wedge \tau_R^n}\|u_n(s)\|^2+\int_0^{t\wedge \tau_R^n}(\|\nabla u_n(s)\|^2+\|u_n(s)\|_{2\sigma+2}^{2\sigma+2})ds\right)\\&\leq& C(\mathbb{E}\|u_n(0)\|^2+1),
\end{eqnarray*}
which completes the proof of Lemma 3.1.\\
{\bf Lemma 3.2} Under the same assumptions as in Theorem 2.5, there exists a constant $C$ such that
\begin{eqnarray}&&\sup\limits_{n}\mathbb{E}\left(\sup\limits_{0\leq t\leq T}\|\nabla u_n(t)\|^2+\int_0^T\|\triangle u_n(s)\|^2ds+
\int_0^T\int_D|u_n(s,x)|^{2\sigma}|\nabla u_n(s,x)|^2dxds\right)\nonumber\\&\leq&C\,(\mathbb{E}\|\nabla u_n(0)\|^2+1).\nonumber
\end{eqnarray}
$Proof:$ We apply It\^{o} formula to $\|\nabla u_n(t)\|^2$ and take the real part of both sides:
\begin{eqnarray}\label{8}\|\nabla u_n(t)\|^2&=&\|\nabla u_n(0)\|^2+2\,\textmd{Re}\int_0^t(\nabla u_n(s),\nabla P_nG(u_n(s)))ds\nonumber\\&&+2\,\textmd{Re}\int_0^t(\nabla u_n(s-),\nabla P_ng(s,u_n(s),z))\tilde{\eta}(ds,dz)\nonumber\\&&
+\int_0^t\int_Z\|\nabla P_ng(s,u_n(s),z)\|^2\eta(ds,dz)\nonumber\\&=&\|\nabla u_n(0)\|^2+I_4+I_5+I_6,
\end{eqnarray}
where
\begin{eqnarray}I_5&=&-2\,\textmd{Re}\int_0^t(\triangle u_n(s), P_ng(s,u_n(s),z)\tilde{\eta}(ds,dz)\nonumber,\\
I_6&=&\int_0^t\int_Z\|\nabla P_ng(s,u_n(s),z)\|^2\eta(ds,dz)\nonumber
\end{eqnarray}
and
\begin{eqnarray}I_4&=&-2\,\textmd{Re}\int_0^t(\triangle u_n(s),(1+i\alpha)\triangle u_n(s)-(1-i\beta)|u_n(s)|^{2\,\sigma}u_n(s) +\gamma u_n(s)+P_n F(u_n(s)))ds\nonumber\\
&=&-2\,\int_0^t\|\triangle u_n(s)\|^2ds+2\,\gamma\int_0^t\|\nabla u_n(s)\|^2ds-2\,\textmd{Re}\int_0^t(\triangle u_n(s),P_n F(u_n(s)))ds\nonumber\\&&+2\,\textmd{Re}(1+i\beta)\int_0^t\int_D|u_n(s)|^{2\sigma}\bar{u}_n(s)\triangle u_n(s) dxds\nonumber\\
&=&-2\,\int_0^t\|\triangle u_n(s)\|^2ds+2\,\gamma\int_0^t\|\nabla u_n(s)\|^2ds+I_{41}+I_{42}.\nonumber
\end{eqnarray}
Take the supremum and expectation in turn over the interval $[0,t]$ on both sides of (\ref{8}), and then we estimate the last three items in the following:
\begin{eqnarray}\label{10}\mathbb{E}\sup\limits_{0\leq s\leq t}I_5(s)&\leq&2\,\mathbb{E} \left[\int_0^t\int_Z\|\nabla u_n(s)\|^2\|\nabla
P_ng(s,u_n(s),z)\|^2\nu(dz)ds\right]^\frac{1}{2}\nonumber
\\
&\leq&2\,\mathbb{E} \left[\int_0^t\int_Z\|\nabla u_n(s)\|^2\int_D|\partial_{u}P_ng(s,u_n(s,x),z)|^2|\nabla u_n(s,x)|^2dx\nu(dz)ds\right]^\frac{1}{2}\nonumber\\
&\leq&2\,\mathbb{E} \left[\int_0^t\int_Z|h(z)|^2\nu(dz)\|\nabla u_n(s)\|^2\int_D|u_n(s,x)|^2|\nabla u_n(s,x)|^2dxds\right]^\frac{1}{2}\nonumber\\
&\leq&2\|h(z)\|_{L^2(Z,\nu)}^2\,\mathbb{E} \left[\int_0^t\|\nabla u_n(s)\|^2\int_D|u_n(s,x)|^2|\nabla u_n(s,x)|^2dxds\right]^\frac{1}{2}\nonumber\\
&\leq&\frac{1}{2}\,\mathbb{E}\sup\limits_{0\leq s\leq t}\|\nabla u_n(s)\|^2+8\|h(z)\|_{L^2(Z,\nu)}^4\,\mathbb{E}\int_0^{t}\int_D|u_n(s,x)|^2|\nabla u_n(s,x)|^2dxds\nonumber\\
&\leq&\frac{1}{2}\,\mathbb{E}\sup\limits_{0\leq s\leq t}\|\nabla u_n(s)\|^2+\epsilon_1\,\mathbb{E}\int_0^t\int_D|u_n(s)|^{2\sigma}|\nabla u_n(s,x)|^2dxds\nonumber\\
&&+C(\epsilon_1, \|h\|_{L^2(Z,\nu)})\,\mathbb{E}\int_0^t\|\nabla u_n(s)\|^2ds.
\end{eqnarray}
The last inequality can be obtained as follows,
\begin{eqnarray}
&&\mathbb{E}\int_0^{t}\int_D|u_n(s,x)|^2|\nabla u_n(s,x)|^2dxds\nonumber\\
&=&\mathbb{E}\int_0^t\int_D|u_n(s,x)|^2|\nabla u_n(s,x)|^{\frac{2}{\sigma}}|\nabla u_n(s,x)|^{2-\frac{2}{\sigma}}dxds\nonumber\\
&\leq&\mathbb{E}\int_0^t\left[\int_D(|u_n(s,x)|^2|\nabla u_n(s,x)|^{\frac{2}{\sigma}})^\sigma dx\right]^\frac{1}{\sigma}\left[\int_D(|\nabla u_n(s,x)|^{2-\frac{2}{\sigma}})^{\frac{\sigma}{\sigma-1}}dx\right]^\frac{\sigma-1}{\sigma}ds\nonumber\\
&\leq&\epsilon\mathbb{E}\int_0^t\int_D|u_n(s,x)|^{2\sigma}|\nabla u_n(s,x)|^2dxds+C(\epsilon)\mathbb{E}\int_0^t\|\nabla u_n(s)\|^2ds.\nonumber
\end{eqnarray}
For the item $I_6$, we notice that
\begin{eqnarray}\label{15}\mathbb{E}\sup\limits_{0\leq s\leq t}I_6(s)&=&\mathbb{E}\int_0^t\int_Z\|\nabla P_ng(s,u_n(s),z)\|^2\eta(ds,dz)\nonumber\\&=&\mathbb{E}\,\int_0^t\int_Z\int_D|\partial_{u}P_ng(s,u_n(s,x),z)|^2|\nabla u_n(s,x)|^2dx\nu(dz)ds,\nonumber\\
&\leq& \mathbb{E}\,\int_0^t\int_D\left[\int_z|h(z)|^2\nu(dz)\right]|u_n(s,x)|^2|\nabla u_n(s,x)|^2dxds\nonumber\\
 &\leq&\|h\|_{L^2(Z,\nu)}^2\mathbb{E}\int_0^t\int_D|u_n(s,x)|^2|\nabla u_n(s,x)|^2dxds\nonumber
 \\&\leq& \epsilon_2\mathbb{E}\int_0^t\int_D|u_n(s,x)|^{2\sigma}|\nabla u_n(s,x)|^2dxds\nonumber\\
&&+C(\epsilon_2, \|h\|_{L^2(Z,\nu)})\mathbb{E}\int_0^t\|\nabla u_n(s,x)\|^2ds.
\end{eqnarray}
 From (\ref{9}), we can obtain
 \begin{eqnarray}&&\mathbb{E}\sup\limits_{0\leq s\leq t}I_{41}(s)\nonumber\\&\leq&2\,\mathbb{E}\left|\int_0^t\int_D\left[\bar{\lambda}_1\cdot\nabla(|u_n(s,x)|^2\bar{u}_n(s,x))+(\bar{\lambda}_2\cdot\nabla \bar{ u}_n(s,x))|u_n(s,x)|^2\right]\triangle u_n(s,x) dxds\right|\nonumber\\
 &\leq&(6|\lambda_1|+2|\lambda_2|)\,\mathbb{E} \int_0^t\int_D|u_n(s,x)|^2|\nabla u_n(s,x)||\triangle u_n(s,x)|dxds\nonumber\\
 &\leq&\frac{1}{2}\,\mathbb{E}\int_0^t\|\triangle u_n(s)\|^2ds+2(3|\lambda_1|+|\lambda_2|)^2\,\mathbb{E}\int_0^t\int_D|u_n(s,x)|^4|\nabla u_n(s,x)|^2dxds\nonumber\\
 &=&\frac{1}{2}\,\mathbb{E}\int_0^t\|\triangle u_n(s)\|^2ds+C(|\lambda_1|,|\lambda_2|)\,\mathbb{E}\int_0^t\int_D|u_n(s,x)|^4|\nabla u_n(s,x)|^\frac{4}{\sigma}|\nabla u_n(s,x)|^{2-\frac{4}{\sigma}}dxds\nonumber\\
 &\leq&\frac{1}{2}\,\mathbb{E}\int_0^t\|\triangle u_n(s)\|^2ds\nonumber\\&&+C\,\mathbb{E}\int_0^t\left[\int_D(|u_n(s,x)|^4|\nabla u_n(s,x)|^{\frac{4}{\sigma}})^\frac{\sigma}{2} dx\right]^\frac{2}{\sigma}\left[\int_D(|\nabla u_n(s,x)|^{2-\frac{4}{\sigma}})^{\frac{\sigma}{\sigma-2}}dx\right]^\frac{\sigma-2}{\sigma}ds\nonumber\\
 &\leq&\frac{1}{2}\,\mathbb{E}\int_0^t\|\triangle u_n(s)\|^2ds\nonumber\\&&+\epsilon_3\,\mathbb{E}\int_0^t\int_D|u_n(s,x)|^{2\sigma}|\nabla u_n(s,x)|^2dx+C(\epsilon_3,|\lambda_1|,|\lambda_2|)\,\mathbb{E}\int_0^t\|\nabla u_n(s)\|^2ds.
 \end{eqnarray}
 In order to estimate the term $I_{42}$, firstly we have
 \begin{eqnarray*} &&2\,\textmd{Re}(1+i\beta)\int_D|u_n(s,x)|^{2\sigma}\bar{u}_n(s,x)\triangle u_n(s,x) dx\nonumber\\
 &=&-2\,\textmd{Re}(1+i\beta)\int_D|u_n(s,x)|^{2(\sigma-1)}[(\sigma+1)|u_n(s,x)|^2|\nabla u_n(s,x)|^2\nonumber\\
 &&\ \ \ \ \ \ \ \ \ \ \ \ \ \ \ \ \ \ \ \ \ \ \ +\sigma (\bar{u}_n(s,x))^2(\nabla u_n(s,x))^2]dx\\
 &=&-\int_D|u_n(s,x)|^{2(\sigma-1)}\sum\limits_{j=1}^2(\bar{u}_n(s,x)\partial_j u_n(s,x),u_n(s,x)\partial_j\bar {u}_n(s,x))\nonumber\\
 &&\ \ \ \ \ \ \ \ \ \ \ \ \ \ \ \ \ \ \ \ \ \ \ \ \ \ \ \ \ \ \ \ \times M(\beta,\sigma)\left(\begin{array}{c}u_n\partial_j\bar{u}_n\\ \bar{u}_n\partial_j u_n \end{array}\right)(s,x)dx,
\end{eqnarray*}
where $$M(\beta,\sigma)=\left(\begin{array}{cc}\sigma+1&(1-i\beta)\sigma\\(1-i\beta)\sigma & \sigma+1\end{array}\right).$$
When $|\beta|<\frac{\sqrt{2\sigma+1}}{\sigma}$, the matrix $M(\beta, \sigma)$ is definitely positive and thus the small eigenvalue $\lambda_\beta$ is positive. This gives
\begin{eqnarray*}2\,\textmd{Re}(1+i\beta)\int_D|u_n(s,x)|^{2\sigma}\bar{u}_n(s,x)\triangle u_n(s,x) dx+2\,\lambda_\beta\int_D|u_n(s,x)|^{2\sigma}|\nabla u_n(s,x)|^2dx\leq 0.
\end{eqnarray*}
Furthermore,
\begin{eqnarray}\label{11}\mathbb{E}\sup\limits_{0\leq s\leq t}I_{42}(s)&\leq&-2\,\lambda_\beta\int_0^t\int_D|u_n(s,x)|^{2\sigma}|\nabla u_n(s,x)|^2dxds.
\end{eqnarray}
Putting (\ref{10})-(\ref{11}) into (\ref{8}), we have
\begin{eqnarray}&&\frac{1}{2}\,\mathbb{E}\sup\limits_{0\leq s\leq t}\|\nabla u_n(s)\|^2+\frac{3}{2}\,\mathbb{E}\int_0^{t}\|\triangle u_n(s)\|^2ds+\frac{C_5}{2}\, \mathbb{E}\int_0^{t}\int_D|u_n(s,x)|^{2\sigma}|\nabla u_n(s,x)|^2dxds\nonumber\\&\leq &\mathbb{E}\|\nabla u_n(0)\|^2+\frac{C_6}{2}\,\mathbb{E}\int_0^{t}\|\nabla u_n(s)\|^2ds,\nonumber\end{eqnarray}
where $\frac{C_5}{2}=2\lambda_\beta-\epsilon_1-\epsilon_2-\epsilon_3$, $\frac{C_6}{2}=C(\epsilon_1,\,\epsilon_2,\,\epsilon_3,\,|\lambda_1|,\,|\lambda_2|,\,\|h\|_{L^2(Z,\nu)})+2\gamma$. Thus, choosing sufficiently small $\epsilon_1,\,\epsilon_2,\,\epsilon_3$ such that $2\lambda_\beta-\epsilon_1-\epsilon_2-\epsilon_3>0$, we can apply Gronwall lemma to get
\begin{eqnarray*}
&&\mathbb{E}\left(\sup\limits_{0\leq s\leq T}\|\nabla u_n(s)\|^2+\int_0^{T}\|\triangle u_n(s)\|^2ds+\int_0^{T}\int_D|u_n(s,x)|^{2\sigma}|\nabla u_n(s,x)|^2dxds\right)\\&\leq&C\,(\mathbb{E}\|\nabla u_n(0)\|^2+1).
\end{eqnarray*}
{\bf Lemma 3.3} If $2\leq p<2\sigma$, there exists a constant $C$ such that
\begin{eqnarray}&&\sup\limits_n\mathbb{E}(\sup\limits_{0\leq s\leq t}\|\nabla u_n(s)\|^p+\int_0^t\|\nabla u_n(s)\|^{p-2}\|\triangle u_n(s)\|^2ds\nonumber\\
&&+ \int_0^t\int_D\|\nabla u_n(s)\|^{p-2}|u_n(s,x)|^{2\sigma}|\nabla u_n(s,x)|^2dxds)\nonumber\\&\leq &C\,\mathbb{E}(\|\nabla u_n(0)\|^p+1).\nonumber\end{eqnarray}
$Proof:$ We apply It\^{o} formula to $\|\nabla u_n(t)\|^p$ and take the real part on both sides to have
\begin{eqnarray}\label{12}\|\nabla u_n(t)\|^p&=&\|\nabla u_n(0)\|^p+p\,\textmd{Re}\int_0^t\|\nabla u_n(s)\|^{p-2}(\nabla u_n(s),\, \nabla P_n G(u_n(s))) ds\nonumber\\&&
+p\,\textmd{Re}\int_0^t\int_Z\|\nabla u_n(s)\|^{p-2}(\nabla u_n(s-),\nabla P_ng(s,u_n(s),z))\tilde{\eta}(ds,dz)\nonumber\\&&+\,\textmd{Re}\int_0^t\int_Z[\,\|\nabla(u_n(s-)+P_ng(s-,u_n(s-),z))\|^p-\|\nabla u_n(s-)\|^p\nonumber\\&&-p\,\|\nabla u_n(s-)\|^{p-2}(\nabla u_n(s-),\nabla P_ng(s-,u_n(s-),z))\,]\eta(ds,dz)\nonumber\\
&=&\|\nabla u_n(0)\|^p+I_7+I_8+I_9.
\end{eqnarray}
Then we estimate the last three items. For $I_7$, we have
\begin{eqnarray}\label{13}I_7(s)&=&-p\int_0^t\|\nabla u_n(s)\|^{p-2}(\triangle u_n(s),\, P_n G(u_n(s))) ds\nonumber\\&=&-p\int_0^t\|\nabla u_n(s)\|^{p-2}\|\triangle u_n(s)\|^2ds+p\,\gamma\int_0^t\|\nabla u_n(s)\|^{p}ds\nonumber\\&&+p\,\textmd{Re}(1+i\beta)\int_0^t\int_D\|\nabla u_n(s)\|^{p-2}|u_n(s,x)|^{2\sigma}\bar{u}_n(s,x)\triangle u_n(s,x)dxds\nonumber\\&&-p\,\textmd{Re}\int_0^t\|\nabla u_n(s)\|^{p-2}(\triangle u_n(s),\, P_n F(u_n(s)))ds\\&=&-p\int_0^t\|\nabla u_n(s)\|^{p-2}\|\triangle u_n(s)\|^2ds+p\,\gamma\int_0^t\|\nabla u_n(s)\|^{p}ds+I_{71}+I_{72}.\nonumber
\end{eqnarray}
A similar estimation to the term $I_{41}$ and $I_{42}$ leads to
\begin{eqnarray}I_{71}(s)&\leq&-2\lambda_\beta\int_0^t\int_D\|\nabla u_n(s)\|^{p-2}|u_n(s,x)|^{2\sigma}|\nabla u_n(s,x)|^2dxds
\end{eqnarray}
and
\begin{eqnarray}\label{14}I_{72}(s)&\leq&\frac{p}{2}\int_0^t\int_D\|\nabla u_n(s)\|^{p-2}\|\triangle u_n(s)\|^2ds+C(\epsilon_4,\,|\lambda_1|,\,|\lambda_2|)\int_0^t\|\nabla u_n(s)\|^pds\nonumber\\&&+\epsilon_4\int_0^t\int_D\|\nabla u_n(s)\|^{p-2}|u_n(s,x)|^{2\sigma}|\nabla u_n(s,x)|^2dxds.
\end{eqnarray}
Taking the supremum and expectation over the interval $[0, t]$ on both sides of (\ref{12}), by (\ref{13})-(\ref{14}) we obtain
\begin{eqnarray}\label{18}\mathbb{E}\sup\limits_{0\leq s\leq t}I_7(s)&\leq&-\frac{p}{2}\,\mathbb{E}\int_0^t\int_D\|\nabla u_n(s)\|^{p-2}\|\triangle u_n(s)\|^2ds+C_7\int_0^t\|\nabla u_n(s)\|^pds\\&&+(-2\lambda_\beta+\epsilon_4)\,\mathbb{E}\int_0^t\int_D\|\nabla u_n(s)\|^{p-2}|u_n(s,x)|^{2\sigma}|\nabla u_n(s,x)|^2dxds,\nonumber
\end{eqnarray}
where $C_7=C(\epsilon_4,\,|\lambda_1|,\,|\lambda_2|)+p\,\gamma$.
Next we will estimate the term $I_8$. Applying the Burkholder-Davis-Gundy inequality and the condition $(C_3)$, we get
\begin{eqnarray}\mathbb{E}\sup\limits_{0\leq s\leq t}I_8(s)&=& p\,\mathbb{E} \sup\limits_{0\leq s\leq t}\left|\int_0^s\int_Z\|\nabla u_n(r)\|^{p-2}(\nabla u_n(r-),\nabla P_ng(r,u_n(r),z))\tilde{\eta}(dr,dz)\right|\nonumber\\&\leq&
\bar{C}\,{\mathbb{E}}\left[\int_0^{t}\int_Z\|\nabla u_n(s)\|^{2p-2}\|\nabla P_ng(s,u_n(s),z)\|^2\nu(dz)ds\right]^\frac{1}{2}\nonumber\\
&\leq&\bar{C}\|h(z)\|_{L^2(Z,\nu)}\,\mathbb{E} \left[\int_0^t\int_D\|\nabla u_n(s)\|^{2p-2}| u_n(s,x)|^2|\nabla u_n(s,x)|^2dxds\right]^\frac{1}{2}\nonumber\\
&\leq&\bar{C}\|h(z)\|_{L^2(Z,\nu)}\nonumber\\
 &&\ \times\mathbb{E} \left[\int_0^t\int_D\|\nabla u_n(s)\|^{p}\|\nabla u_n(s)\|^{p-2}|u_n(s,x)|^2|\nabla u_n(s,x)|^2dxds\right]^\frac{1}{2}\nonumber\\&\leq&
\frac{1}{2}\,\mathbb{E}\sup\limits_{0\leq s\leq t}\|\nabla u_n(s)\|^p+C\,\mathbb{E}\int_0^{t}\int_D\|\nabla u_n(s)\|^{p-2}|u_n(s,x)|^2|\nabla u_n(s,x)|^2dxds\nonumber\\&\leq&\frac{1}{2}\,\mathbb{E}\sup\limits_{0\leq s\leq t}\|\nabla u_n(s)\|^p+C(\epsilon_5)\|\nabla u_n(s)\|^2)ds\nonumber\\&&+C\mathbb{E}\int_0^{t}\|\nabla u_n(s)\|^{p-2}(\epsilon_5\int_D|u_n(s,x)|^{2\sigma}|\nabla u_n(s,x)|^2dx\nonumber\\&\leq&\frac{1}{2}\,\mathbb{E}\sup\limits_{0\leq s\leq t}\|\nabla u_n(s)\|^p+C(\epsilon_5,\,\|h(z)\|_{L^2(Z,\nu)})\,\mathbb{E}\int_0^t\|\nabla u_n(s)\|^pds\nonumber\\&&+\epsilon_5\,\mathbb{E}\int_0^{t}\int_D\|\nabla u_n(s)\|^{p-2}|u_n(s,x)|^{2\sigma}|\nabla u_n(s,x)|^2dxds.
\end{eqnarray}
On the other hand, by the Taylor formula, we have
\begin{eqnarray*}
|\|x+h\|^p-\|x\|^p-p\|x\|^{p-2}(x,h)|\leq C_p(\|x\|^{p-2}\|h\|^2+\|h\|^p)\ \ \ \ \ \ {\rm for\ all}\ x,h\in H_n.
\end{eqnarray*}
Thus
\begin{eqnarray}\label{16}&&\mathbb{E}\sup\limits_{0\leq s\leq t}I_9(s)\nonumber\\
&\leq&C\mathbb{E}\int_0^t\int_Z(\|\nabla u_n(s)\|^{p-2}\|\nabla P_ng(s,u_n(s),z)\|^2\nu(dz)ds\nonumber\\
&&+C\mathbb{E}\int_0^t\int_Z\|\nabla P_ng(s,u_n(s),z)\|^p)\nu(dz)ds\nonumber\\&=&I_{91}(t)+I_{92}(t).
\end{eqnarray}
From the condition $(C_3)$ and (\ref{15}), we deduce
\begin{eqnarray}I_{91}(t)&\leq&C(\|h\|_{L^2(Z,\nu)}^2)\,\mathbb{E}\int_0^t\int_D\|\nabla u_n(s)\|^{p-2}|u_n(s,x)|^2|\nabla u_n(s,x)|^2dxds\nonumber\\&\leq&\epsilon_6\,\mathbb{E}\int_0^{t}\int_D\|\nabla u_n(s)\|^{p-2}|u_n(s,x)|^{2\sigma}|\nabla u_n(s,x)|^2dxds\nonumber\\&&+C(\epsilon_6,\,\|h(z)\|_{L^2(Z,\nu)})\,\mathbb{E}\int_0^t\|\nabla u_n(s)\|^pds.
\end{eqnarray}
Now we estimate the last term of (\ref{16}). By Young inequality, it follows that
\begin{eqnarray}\label{17}I_{92}(t)&=&\|h\|_{L^p(Z,\nu)}^p\,\mathbb{E}\int_0^t\left(\int_D|u_n(s,x)|^2|\nabla u_n(s,x)|^2 dx\right)^{\frac{p}{2}}ds\nonumber\\
&=&C\,\mathbb{E}\int_0^t\left(\int_D|u_n(s,x)|^2|\nabla u_n(s,x)|^{\frac{4}{p}} |\nabla u_n(s,x)|^{2-\frac{4}{p}}dx\right)^{\frac{p}{2}}ds\nonumber\\
&\leq&C\,\mathbb{E}\int_0^t\left(\int_D(|u_n(s,x)|^2|\nabla u_n(s,x)|^{\frac{4}{p}})^{\frac{p}{2}}dx\right)\nonumber\\
&&\ \ \ \ \ \ \ \times\left(\int_D|\nabla u_n(s,x)|^{2-\frac{4}{p}})^{\frac{p}{p-2}}dx\right)^{\frac{p-2}{2}}ds\nonumber\\
&=&C\,\mathbb{E}\int_0^t\|\nabla u_n(s)\|^{p-2}\int_D|u_n(s,x)|^p|\nabla u_n(s,x)|^2dxds\nonumber\\&=&C\,\mathbb{E}\int_0^t\|\nabla u_n(s)\|^{p-2}
\int_D|u_n(s,x)|^p|\nabla u_n(s,x)|^{\frac{p}{\sigma}}|\nabla u_n(s,x)|^{2-\frac{p}{\sigma}}dxds\nonumber\\&\leq&\mathbb{E}\int_0^t\|\nabla u_n(s)\|^{p-2}
\left(\epsilon_7\int_D|u_n(s,x)|^{2\sigma}|\nabla u_n(s,x)|^2dx+C(\epsilon_7)\|\nabla u_n(s)\|^2\right)ds\nonumber\\&=&\epsilon_7\,\mathbb{E}\int_0^{t}\int_D\|\nabla u_n(s)\|^{p-2}|u_n(s,x)|^{2\sigma}|\nabla u_n(s,x)|^2dxds\nonumber\\&&+C(\epsilon_7)\,\mathbb{E}\int_0^t\|\nabla u_n(s)\|^pds.
\end{eqnarray}
Combining (\ref{18})-(\ref{17}) with (\ref{12}), we infer that
\begin{eqnarray}&&\frac{1}{2}\,\mathbb{E}\sup\limits_{0\leq s\leq t}\|\nabla u_n(s)\|^p+\frac{p}{2}\,\mathbb{E}\int_0^t\|\nabla u_n(s)\|^{p-2}\|\triangle u_n(s)\|^2ds\nonumber\\&&+\frac{C_8}{2}\,\mathbb{E}\int_0^t\int_D\|\nabla u_n(s)\|^{p-2}|u_n(s,x)|^{2\sigma}|\nabla u_n(s,x)|^2dxds\nonumber\\&\leq &\mathbb{E}\|u_n(0)\|^2+\frac{C_9}{2}\,\mathbb{E}\int_0^t\|\nabla u_n(s)\|^pds,\nonumber\end{eqnarray}
where $\frac{C_8}{2}=2\lambda_\beta-\sum_{i=4}^7\epsilon_i$, $\frac{C_9}{2}=p\,\gamma +C(|\lambda_1|,\,|\lambda_2|,\,\epsilon_i), i=4,5,6,7$. Thus, choosing sufficiently small $\epsilon_i$ such that $2\lambda_\beta-\sum_{i=4}^7\epsilon_i>0$, we can use Gronwall lemma to get
\begin{eqnarray*}
&&\mathbb{E}(\sup\limits_{0\leq s\leq T}\|\nabla u_n(s)\|^p+\int_0^{T}\int_D\|\nabla u_n(s)\|^{p-2}\|\triangle u_n(s)\|^2ds\\&& +\int_0^{T}\|\nabla u_n(s)\|^{p-2}|u_n(s,x)|^{2\sigma}|\nabla u_n(s,x)|^2dxds)\\&\leq&C\,(\mathbb{E}\|\nabla u_n(0)\|^p+1).
\end{eqnarray*}
Since the constant $C$ is independent of $n$, the proof of Lemma 3.3 is complete.\\
\indent Step 2 :  We study the weak convergence of approximating sequences in this step. Under the same assumptions as in Theorem 2.5, there exists a subsequence of $\{u_n\}$, still denoted by $\{u_n\}$, satisfying the following Lemma:\\
{\bf Lemma 3.4 } Let $T(u)=-(1-i\beta)|u|^{2\sigma}u$, $s=(2\sigma+2)^*=\frac{2\sigma+2}{2\sigma+1}$, $2\leq p<2\sigma$. There exist processes $\tilde{u}\in L^2([0,T]\times\Omega;H^2)\cap L^{p}(\Omega;L^{\infty}([0,T];V))\cap L^{2\sigma+2}([0,T]\times\Omega;L^{2\sigma+2}(D))$, $\tilde{T}\in L^s([0,T]\times\Omega;L^s(D))$, $\tilde{F}\in L^2([0,T]\times\Omega;H^{-2})$, $Y\in L^2([0,T]\times\Omega;L^2(Z,\nu;H))$, where $H^{-2}$ is the dual space of $H^2$, such that\\
 \indent (1) $u_n \rightharpoonup \tilde{u}$ \,in \,$L^2([0,T]\times\Omega;H^2)\cap L^{2\sigma+2}([0,T]\times\Omega;L^{2\sigma+2}(D))$;\\
 \indent (2) $u_n$ is weak star converging to $\tilde{u}$ in $L^{p}(\Omega;L^{\infty}([0,T];V))$;\\
 \indent (3) $P_n T(u_n)\rightharpoonup \tilde{T}$\,in \,$L^s([0,T]\times\Omega;L^s(D))$;\\
 \indent (4) $P_n F(u_n)\rightharpoonup \tilde{F}$ in $L^2([0,T]\times\Omega;H^{-2})$;\\
 \indent (5) $P_ng(s,u_n(s),z)\rightharpoonup Y$ in $L^2([0,T]\times\Omega;L^2(Z,\nu;H))$.\\
 $Proof:$ It is easy to see that (1) and (3) are straightforward consequences of Lemma 3.1.
Since $ L^{p}(\Omega; L^{\infty}([0,T]; V))\cong(L^{p^{*}}(\Omega; L^1([0,T]; V^{*})))^{*}$, according to the Banach-Alaoglu theorem, we can get another weakly star convergent subsequence, still denoted by $\{u_n\}$, and $ \tilde{u}\in  L^{p}(\Omega; L^{\infty}([0,T]; V))$ such that (2) holds. Moreover, for $v\in L^{2}\,([0,T]\times\Omega;H^2)$, we find that
 \begin{eqnarray}
&&\mathbb{E}\int_0^T\|P_n F(u_n(s))\|_{H^{-2}}^2ds\nonumber\\
&\leq&3(|\lambda_1|+|\lambda_2|)\,\sup\limits_{\|v\|_{L^2([0,T]\times\Omega;H^2)}\leq 1}\mathbb{E}\int_0^T\int_D| u_n(s,x)|^2|\nabla u_n(s,x)||v(s,x)|dxds\nonumber\\
&\leq&3(|\lambda_1|+|\lambda_2|)\,\sup\limits_{\|v\|_{L^2([0,T]\times\Omega;H^2)}\leq 1}\mathbb{E}\int_0^T\int_D| u_n(s,x)|^2|\nabla u_n(s,x)|^{\frac{2}{\sigma}}|\nabla u_n(s,x)|^{1-\frac{2}{\sigma}}\nonumber\\
&&\ \ \ \ \ \ \ \ \ \ \ \ \ \ \ \ \ \ \ \ \ \ \ \ \ \ \ \ \ \ \ \ \ \ \ \ \ \ \ \ \ \ \ \ \ \ \ \ \times|v(s,x)|dxds\nonumber\\
&\leq&3(|\lambda_1|+|\lambda_2|)\,\sup\limits_{\|v\|_{L^2([0,T]\times\Omega;H^2)}\leq 1}\mathbb{E}\int_0^T\left(\int_D| u_n(s,x)|^{2\sigma}|\nabla u_n(s,x)|^2dx\right)^{\frac{1}{\sigma}}\nonumber\\
&&\ \ \ \ \ \ \ \ \ \ \ \ \ \ \ \ \ \ \ \ \ \ \ \ \ \ \ \ \ \ \ \ \ \ \ \ \ \ \ \ \ \ \ \ \ \ \times\left(\int_D|\nabla u_n(s,x)|^2dx\right)^{\frac{\sigma-2}{2\sigma}}\|v(s)\|ds\nonumber\\
&\leq& \sup\limits_{\|v\|_{L^2([0,T]\times\Omega;H^2)}\leq 1}\mathbb{E}\int_0^T\|v\|^2ds+\mathbb{E}\int_0^T\int_D| u_n(s,x)|^{2\sigma}|\nabla u_n(s,x)|^2dx\nonumber\\
&&+C(|\lambda_1|,\,|\lambda_2|)\mathbb{E}\int_0^T\|\nabla u_n(s)\|^2ds<\infty,\nonumber
\end{eqnarray}
which implies that (4) holds. As for (5), by  the condition $(C_1)$ and Lemma 3.1, we have
$$ \mathbb{E}\int_0^T\int_Z\|P_ng(s,u_n(s),z)\|^2\nu(dz)ds\leq\mathbb{E}\int_0^T(k_1\|u_n(s)\|^2+k_2\,\|\nabla u_n(s)\|^2)ds<\infty.$$
So we konw that $P_ng$ has a subsequence converging weakly in $L^2([0,T]\times\Omega, L^2(Z,\nu;H))$ to $Y$.\\
\indent Step 3 : Define $X:=(1+i\alpha)\triangle\tilde{u}+\tilde{T}+\gamma \tilde{u}+\tilde{F}$ and then a $V^{*}$-valued process $u$ by
$$u(t):=u_0+\int_0^t X(s)ds+\int_0^t\int_Z Y(s,z)\tilde{\eta}(dz,ds)$$
such that $u$ is a $V^{*}$-valued  modification of the $V$-valued process $\tilde{u}$ in Lemma 3.4, i.e $u=\tilde{u}$ $dt\times \mathbb{P}$-a.e. in $V$. The aim of this step is to verify the following identities $dt\times \mathbb{P}$-a.e. on $[0,T]\times\Omega$:\\
$$G(\tilde{u}(s))=X(s),\ \ \ \ g(s,\tilde{u},z)=Y(s,z)\ \ \ \ {\rm a.s.}$$

Take $\phi$ from $L^2([0,T]\times\Omega;V)\cap L^{2\sigma+2}([0,T]\times\Omega;L^{2\sigma+2}(D))\cap L^p (\Omega;L^{\infty}([0,T]; V))$. Denote  \begin{eqnarray*}&&I=-(1-i\beta)\langle\,|u_n|^{2\sigma}u_n-|\phi|^{2\sigma}\phi,u_n-\phi\,\rangle,\\
&&m=\textmd{Re}\langle\,|u_n|^{2\sigma}u_n-|\phi|^{2\sigma}\phi,u_n-\phi\,\rangle,\,\,\ \ \ n
=\textmd{Im}\langle\,|u_n|^{2\sigma}u_n-|\phi|^{2\sigma}\phi,u_n-\phi\,\rangle,\\
&&F(u_n)=((2\lambda_1+\lambda_2)\cdot\nabla u_n)|u_n|^2+(\lambda_1\cdot\nabla u_n)u_n^2=F_1(u_n)+F_2(u_n),\\
&&J=(F_1(u_n)-F_1(\phi),u_n-\phi)=((2\lambda_1+\lambda_2)\cdot(|u_n|^2\nabla u_n-|\phi|^2\nabla \phi),u_n-\phi),\\
&&K=(F_2(u_n)-F_2(\phi),u_n-\phi)=(\lambda_1\cdot(u_n^2\nabla u_n-\phi^2\nabla \phi),u_n-\phi).\end{eqnarray*}
In order to obtain monotonicity properties of $G$, we need the following two lemmas.\\
 {\bf Lemma 3.5 }If $0<|\beta|<\frac{\sqrt{2\sigma+1}}{\sigma}$, we have
 $$\textmd{Re }I\leq-(1-\frac{\sigma}{\sqrt{2\sigma+1}}|\beta|)2^{-2\sigma}\|u_n-\phi\|_{2\sigma+2}^{2\sigma+2}.$$
 $Proof:$  Note that\begin{eqnarray*} m&=&\textmd{Re}\langle\,|u_n|^{2\sigma}u_n-|\phi|^{2\sigma}\phi,u_n-\phi\,\rangle\\
 &=&\textmd{Re}\langle\,(|u_n|^{2\sigma}-|\phi|^{2\sigma})u_n+|\phi|^{2\sigma}(u_n-\phi), u_n-\phi\,\rangle\\&=&\int_D|\phi(s,x)|^{2\sigma}|u_n(s,x)-\phi(s,x)|^2dx\nonumber\\
&&+\frac{1}{2}\int_D(|u_n(s,x)|^{2\sigma}-|\phi(s,x)|^{2\sigma})(|u_n(s,x)|^2+|u_n(s,x)-\phi(s,x)|^2-|\phi(s,x)|^2)dx\\
 &\geq&\int_D|\phi(s,x)|^{2\sigma}|u_n(s,x)-\phi(s,x)|^2dx\nonumber\\
&&+\frac{1}{2}\int_D(|u_n|^{2\sigma}-|\phi(s,x)|^{2\sigma})|u_n(s,x)-\phi(s,x)|^2dx\\&=&\frac{1}{2}\int_D(|u_n(s,x)|^{2\sigma}+|\phi(s,x)|^{2\sigma})|u_n(s,x)-\phi(s,x)|^2dx\\
 &\geq&2^{-2\sigma}\,\int_D|u_n(s,x)-\phi(s,x)|^{2\sigma+2}dx.
 \end{eqnarray*}
Let $p=2\sigma+2$, then $\frac{|n|}{m}\le \frac{\sigma}{\sqrt{2\sigma+1}}$. By Lemma 2.3, we find that
\begin{eqnarray*}\textmd{Re}\,I&=&-m-\beta\,n\leq -(1-\frac{\sigma}{\sqrt{2\sigma+1}}|\beta|)\,m\\&\leq&
-(1-\frac{\sigma}{\sqrt{2\sigma+1}}|\beta|)2^{-2\sigma}\|u_n-\phi\|_{2\sigma+2}^{2\sigma+2},
 \end{eqnarray*}
if \,$0<|\beta|<\frac{\sqrt{2\sigma+1}}{\sigma}$.\\
 {\bf Lemma 3.6} We denote $w=u_n-\phi$ for simplicity, there exist enough small parameters $\tilde{\epsilon},\,\hat{\epsilon}$
 such that
  \begin{eqnarray}\textmd{Re}\,J&\leq&\tilde{\epsilon}\|\nabla w\|_2^2+\hat{\epsilon}\|w\|_{2\sigma+2}^{2\sigma+2}+\big(C(\epsilon_8,\epsilon_9)+(\epsilon_{13}+\epsilon_{15})\|\phi\|^2+C(\epsilon_{10},\,\epsilon_{11})\|\nabla \phi\|^{\frac{2\sigma}{\sigma-1}}\nonumber\\&&+C(\epsilon_{12},\epsilon_{13})\|\nabla\phi\|^{\frac{7\sigma-2}{\sigma+1}}+C(\epsilon_{14},\epsilon_{15})\|\nabla\phi\|^{\frac{10\sigma+4}{\sigma+4}}\big)\|w\|^2.\nonumber\end{eqnarray}
$Proof:$ First note that
\begin{eqnarray}\label{21}\textmd{Re}\,J&\leq&C\int_D|(|u_n(s,x)|^2\nabla u_n(s,x)-|\phi(s,x)|^2\nabla \phi(s,x))\bar{w}(s,x)|dx\nonumber\\&=&C\int_D|(|u_n(s,x)|^2\nabla w(s,x)+u_n(s,x)\nabla\phi(s,x)\bar{w}(s,x)\nonumber\\
&&\ \ \ \ \ \ \ \ \ +\bar{\phi}(s,x)\nabla\phi(s,x) w(s,x))\bar{w}(s,x)|dx\nonumber\\&\leq&C\int_D(|w(s,x)|^3|\nabla w(s,x)|+|w(s,x)|^3|\nabla \phi(s,x)|+|w(s,x)|^2
|\phi(s,x)||\nabla\phi(s,x)|\nonumber\\
&&\ \ \ \ \ \ \ \ \ +|w(s,x)||\nabla w(s,x)||\phi(s,x)|^2)dx\nonumber\\&=&J_1+J_2+J_3+J_4.\end{eqnarray}
In the following, we will estimate the four items $J_1,\ J_2,\ J_3,\ J_4$ by Sobolev embedding inequality, H\"{o}lder inequality and Young inequality. For $J_1$, we have
\begin{eqnarray}\label{19}J_1&=&\int_D|w(s,x)|^3|\nabla w(s,x)|dx\leq\|\nabla w\|\|w\|_6^3\leq\epsilon_8\|\nabla w\|^2+C(\epsilon_8)\|w\|_6^6\nonumber\\
&\leq&\epsilon_8\|\nabla w\|^2+\epsilon_9\|w\|_{2\sigma+2}^{2\sigma+2}+C(\epsilon_8,\,\epsilon_9)\|w\|^2.
\end{eqnarray}
The last inequality of (\ref{19}) is due to (\ref{20}). To estimate the second item $J_2$, we first notice
\begin{eqnarray}\|w\|_6\leq\|w\|^{1-\theta}\|w\|_{2\sigma+2}^{\theta},\nonumber
  \end{eqnarray}
where $\theta=\frac{2\sigma+2}{3\sigma}$. From the Gagliardo-Nirenberg inequality for two-dimensional domain, we have
\begin{eqnarray}\|w\|_6\leq\|w\|^{\frac{1}{3}}\|\nabla w\|^{\frac{2}{3}}.\nonumber
 \end{eqnarray}
Thus
\begin{eqnarray}J_2&=&\int_D|w(s,x)|^3|\nabla \phi(s,x)|dx\leq\|w\|_6^{\frac{3}{2}}\|w\|_6^{\frac{3}{2}}\|\nabla\phi\|\nonumber\\&\leq&(\|w\|^{\frac{1}{3}}\|\nabla w\|^{\frac{2}{3}})^{\frac{3}{2}}\|w\|^{\frac{\sigma-2}{2\sigma}}\|w\|_{2\sigma+2}^{\frac{2\sigma+2}{2\sigma}}\|\nabla\phi\|\nonumber\\
&=&\|\nabla w\|\|w\|_{2\sigma+2}^{\frac{2\sigma+2}{2\sigma}}\|w\|^{\frac{\sigma-1}{\sigma}}\|\nabla\phi\|\nonumber\\
&\leq&\epsilon_{10}\|\nabla w\|^2+\epsilon_{11}\|w\|_{2\sigma+2}^{2\sigma+2}+C(\epsilon_{10},\,\epsilon_{11})\|\nabla\phi\|^{\frac{2\sigma}{\sigma-1}}\|w\|^2.
\end{eqnarray}
For the third item $J_3$ in (\ref{21}), we deduce that
\begin{eqnarray}J_3&=&\int_D|w(s,x)|^2|\phi(s,x)||\nabla\phi|dx\leq\|\nabla \phi\|\|\phi\|_{\frac{8\sigma-4}{\sigma-2}}\|w\|_{\frac{16\sigma-8}{3\sigma}}\nonumber\\&\leq&\|\nabla \phi\|\|\phi\|^{\frac{\sigma-2}{4\sigma-2}}\|\nabla \phi\|^{\frac{3\sigma}{4\sigma-2}}\|w\|^{\frac{3\sigma}{4\sigma-2}}\|\nabla w\|^\frac{5\sigma-4}{4\sigma-2}\nonumber\\&\leq&\epsilon_{12}\|\nabla w\|^2+C(\epsilon_{12})\big(\|\nabla \phi\|^{\frac{7\sigma-2}{4\sigma-2}}\| \phi\|^{\frac{\sigma-2}{4\sigma-2}}\big)^{\frac{8\sigma-4}{3\sigma}}\|w\|^2\nonumber\\
&=&\epsilon_{12}\|\nabla w\|^2+C(\epsilon_{12})\|\nabla \phi\|^{\frac{2(7\sigma-2)}{3\sigma}}\| \phi\|^{\frac{2(\sigma-2)}{3\sigma}}\|w\|^2\nonumber\\&\leq&\epsilon_{12}\|\nabla w\|_2^2+\big(\epsilon_{13}\|\phi\|_2^2+C(\epsilon_{12},\,\epsilon_{13} )\|\nabla \phi\|^{\frac{7\sigma-2}{\sigma+1}}\big)\|w\|^2,
\end{eqnarray}
where $2<\frac{7\sigma-2}{\sigma+1}<2\sigma.$ \\
Now we turn to the last item $J_4$. It follows that
\begin{eqnarray}\label{22}J_4&=&\int_D|w(s,x)||\nabla w(s,x)||\phi|^2dx\leq\|\nabla w\|\|w\|_{\frac{3\sigma}{\sigma+1}}\|\phi\|_{\frac{12\sigma}{\sigma-2}}^2\nonumber\\
&\leq&\|\nabla w\|^{1+\frac{\sigma-2}{3\sigma}}\|w\|^{\frac{2\sigma+2}{3\sigma}}\|\phi\|_{\frac{12\sigma}{\sigma-2}}^2\nonumber\\
&\leq&\epsilon_{14}\|\nabla w\|^2+C(\epsilon_{14})\|\phi\|_{\frac{12\sigma}{\sigma-2}}^{\frac{6\sigma}{\sigma+1}}\|w\|^2\nonumber\\
&\leq&\epsilon_{14}\|\nabla w\|^2+C(\epsilon_{14})\big(\|\phi\|^{\frac{\sigma-2}{6\sigma}}\|\nabla \phi\|^{\frac{5\sigma+2}{6\sigma}}\big)^{\frac{6\sigma}{\sigma+1}}\|w\|^2\nonumber\\
&=&\epsilon_{14}\|\nabla w\|^2+C(\epsilon_{14})\|\phi\|^{\frac{\sigma-2}{\sigma+1}}\|\nabla \phi\|^{\frac{5\sigma+2}{\sigma+1}}\|w\|^2\nonumber\\&\leq&\epsilon_{14}\|\nabla w\|^2+\big(\epsilon_{15}\|\phi\|^2+C(\epsilon_{14},\,\epsilon_{15})\|\nabla \phi\|^{\frac{10\sigma+4}{\sigma+4}}\big)\|w\|^2,\end{eqnarray}
where $2<\frac{10\sigma+4}{\sigma+4}<2\sigma.$\\
Based on (\ref{21})--(\ref{22}), we deduce that
\begin{eqnarray}\textmd{Re}\,J&\leq&\tilde{\epsilon}\|\nabla w\|^2+\hat{\epsilon}\|w\|_{2\sigma+2}^{2\sigma+2}+\big(C(\epsilon_8,\epsilon_9)+(\epsilon_{13}+\epsilon_{15})\|\phi\|^2+C(\epsilon_{10},\,\epsilon_{11})\|\nabla \phi\|^{\frac{2\sigma}{\sigma-1}}\nonumber\\&&+
C(\epsilon_{12},\epsilon_{13})\|\nabla\phi\|^{\frac{7\sigma-2}{\sigma+1}}+C(\epsilon_{14},\epsilon_{15})\|\nabla\phi\|^{\frac{10\sigma+4}{\sigma+4}}\big)\|w\|^2,\nonumber\end{eqnarray}
where $\tilde{\epsilon}=\epsilon_{8}+\epsilon_{10}+\epsilon_{12}+\epsilon_{14}$ and $\hat{\epsilon}=\epsilon_{9}+\epsilon_{11}$. This puts an end of the proof of Lemma 3.6.

The estimate of $\textmd{Re}\,K$ is similar to Lemma 3.6. From Lemma 3.5 and Lemma 3.6, we know
\begin{eqnarray}\label{27}&&
\textmd{Re}\,\langle\, P_n G(u_n)-P_n G(\phi),u_n-\phi\,\rangle\nonumber\\&=&\textmd{Re}\,[(1+i\alpha)\langle\,\triangle(u_n-\phi), u_n-\phi\,\rangle+\gamma\langle\,u_n-\phi, u_n-\phi\,\rangle+I+J+K]\nonumber\\&=&(-1+ \tilde{\epsilon})\|\nabla(u_n-\phi)\|^2+\left(-(1-\frac{\sigma}{\sqrt{2\sigma+1}}|\beta|)2^{-2\sigma}+\hat{\epsilon}\right)\|u_n-\phi\|_{2\sigma+2}^{2\sigma+2}\nonumber\\
&&+\big(C(\epsilon_8,\,\epsilon_9,\gamma)+(\epsilon_{13}+\epsilon_{15})\|\phi\|^2+C(\epsilon_{10},\,\epsilon_{11})\|\nabla \phi\|^{\frac{2\sigma}{\sigma-1}}\nonumber\\&&+
C(\epsilon_{12},\epsilon_{13})\|\nabla\phi\|^{\frac{7\sigma-2}{\sigma+1}}+C(\epsilon_{14},\epsilon_{15})\|\nabla\phi\|^{\frac{10\sigma+4}{\sigma+4}}\big)\|(u_n-\phi)\|^2.\nonumber
\end{eqnarray}
Applying It\^{o} formula to the process $e^{-r(t)}\|u_n(t)\|^2$ ($r(t)$ will be defined later) and taking the real part, we obtain
\begin{eqnarray}
e^{-r(t)}\|u_n(t)\|^2&=&\|u_n(0)\|^2-\int_0^t e^{-r(s)}r^{'}(s)\|u_n(s)\|^2ds\nonumber\\&&+2\,\textmd{Re}\int_0^t e^{-r(s)}\langle u_n(s), P_n G(u_n(s))\rangle ds\nonumber\\&&+2\,\textmd{Re}\int_0^t e^{-r(s)}(u_n(s-),P_ng(s,u_n(s),z))\tilde{\eta}(dz,ds)\nonumber\\&&+\int_0^t\int_Ze^{-r(s)}\|P_ng(s,u_n(s),z)\|^2\eta(ds,dz).\nonumber
\end{eqnarray}
Thus by taking the expectation on both sides of above equality, we have
\begin{eqnarray}\label{23}&&\mathbb{E}\left(e^{-r(t)}\|u_n(t)\|^2\right)-\mathbb{E}\|u_n(0)\|^2\nonumber\\&=&\mathbb{E}\Big(\textmd{Re}(-\int_0^t e^{-r(s)}r^{\prime}(s)\left(2(u_n(s),\phi(s))-\|\phi(s)\|^2\right)ds\nonumber\\&&\ \ \ +2\int_0^t e^{-r(s)}\left(\langle\, P_n G(u_n(s))-P_n G(\phi(s)),\phi(s)\,\rangle +\langle\, P_n G(\phi(s)),u_n(s)\,\rangle\right)ds\nonumber\\&&\ \ \ +\int_0^t\int_Ze^{-r(s)}\left(2(P_ng(s,u_n(s),z),P_ng(\phi(s)))-\|P_ng(s,
\phi(s),z)\|^2\right)\nu(dz)ds\nonumber\\
&&\ \ \ -\int_0^t e^{-r(s)}r^{\prime}(s)\|u_n(s)-\phi(s)\|^2ds\nonumber\\&&\ \ \ +2\int_0^t e^{-r(s)}\langle \,P_n G(u_n(s))-P_n G(\phi(s)),u_n(s)-\phi(s)\,\rangle ds\nonumber\\&&\ \ \ +\int_0^t\int_Ze^{-r(s)}\|P_ng(s,u_n(s),z)-P_ng(s,\phi(s),z)\|^2\nu(dz)ds\Big).
\end{eqnarray}
According to the inequality (\ref{23}) and the condition $(C_2)$, we can set
 \begin{eqnarray}\label{26}r(t)&=&\int_0^t[2(C(\epsilon_8,\,\epsilon_9,\gamma)+(\epsilon_{13}+\epsilon_{15})\|\phi(s)\|^2+C(\epsilon_{10},\,\epsilon_{11})\|\nabla \phi(s)\|^{\frac{2\sigma}{\sigma-1}}\nonumber\\&&+C(\epsilon_{12},\epsilon_{13})\|\nabla\phi(s)\|^{\frac{7\sigma-2}{\sigma+1}}+C(\epsilon_{14},\epsilon_{15})\|\nabla\phi(s)\|^{\frac{10\sigma+4}{\sigma+4}})+k_3]ds.\end{eqnarray}
Then we have
\begin{eqnarray}\label{34}&&-r^{\prime}(s)\|u_n(s)-\phi(s)\|^2+2\,\textmd{Re}\,\langle\, P_n G(u_n(s))-P_n G(\phi(s)),u_n(s)-\phi(s)\,\rangle\nonumber\\&&+\int_Z\|P_ng(s,u_n(s),z)-P_ng(s,\phi(s),z)\|^2\nu(dz)\nonumber\\&\leq&(-2+2\tilde{\epsilon}+k_4)\|\nabla(u_n(s)-\phi(s))\|^2+2K(\|\nabla(u_n(s)-\phi(s))\|_{2\sigma+2}^{2\sigma+2}\nonumber\\&\leq&0,
\end{eqnarray}
 provided that $\tilde{\epsilon},\,\hat{\epsilon},\,k_4$ are small enough such that $K=-(1-\frac{\sigma}{\sqrt{2\sigma+1}}|\beta|)2^{-2\sigma}+\hat{\epsilon}< 0$ and $-2+2\tilde{\epsilon}+k_4< 0$. Thus,\\
\begin{eqnarray}&&\mathbb{E}\,\textmd{Re}(-\int_0^t e^{-r(s)}r^{\prime}(s)\|u_n(s)-\phi(s)\|^2ds\nonumber\\
&&+2\int_0^t e^{-r(s)}\langle \,P_n G(u_n(s))-P_n G(\phi(s)),u_n(s)-\phi(s)\,\rangle ds\nonumber\\&&\,\,\,+\int_0^t\int_Ze^{-r(s)}\|P_ng(s,u_n(s),z)-P_ng(s,\phi(s),z)\|^2\nu(dz)ds)\leq 0.\nonumber
\end{eqnarray}
Now we can rewrite (\ref{23}) as below:
 \begin{eqnarray}
 &&\mathbb{E}\left(e^{-r(t)}\|u_n(t)\|^2\right)-\mathbb{E}\|u_n(0)\|^2\nonumber\\&\leq&-\mathbb{E}\,\textmd{Re}\int_0^t e^{-r(s)}r^{\prime}(s)\left(2(u_n(s),\phi(s))-\|\phi(s)\|^2\right)ds\nonumber\\&&+2\,\mathbb{E}\,\textmd{Re}\int_0^t e^{-r(s)}\left(\langle\, P_n G(u_n(s))-P_n G(\phi(s)),\phi(s)\,\rangle +\langle\, P_n G(\phi(s)),u_n(s)\,\rangle\right)ds\nonumber\\&&+\mathbb{E}\,\textmd{Re}\int_0^t\int_Ze^{-r(s)}\left(2(P_ng(s,u_n(s),z),P_ng(s,\phi(s),z))-\|P_ng(
s,\phi(s),z)\|^2\right)\nu(dz)ds.\nonumber
\end{eqnarray}
 By lower semi-continuity property of weak convergence and Lemma 3.5, we have
 \begin{eqnarray}\label{24}&&\mathbb{E}\left(e^{-r(T)}\|u(T)\|^2-\|u_n(0)\|^2\right)\nonumber\\&\leq&\liminf\limits_{n\rightarrow\infty}\mathbb{E}
 \left(e^{-r(T)}\|u_n(T)\|^2-
 \|u_n(0)\|^2\right)\nonumber\\&\leq&\liminf\limits_{n\rightarrow\infty}\mathbb{E}\,\textmd{Re}(\int_0^T e^{-r(s)}r^{\prime}(s)\left(\|\phi\|^2- 2(u_n(s),\phi(s))\right)ds\nonumber\\
 &&+2\int_0^T e^{-r(s)}\left(\langle\, P_n G(u_n(s))-P_n G(\phi(s)),\phi(s)\,\rangle +\langle\, P_n G(\phi(s)),u_n(s)\,\rangle\right)ds\nonumber\\&&+\int_0^T\int_Ze^{-r(s)}\left(2(P_ng(s,u_n(s),z),P_ng(s,\phi(s),z))-\|P_ng(
s,\phi(s),z)\|^2\right)\nu(dz)ds\nonumber\\&=&\mathbb{E}\,\textmd{Re}[\int_0^T e^{-r(s)}r^{\prime}(s)\left(\|\phi(s)\|^2- 2(\tilde{u}(s),\phi(s))\right)ds\nonumber\\&&+2\int_0^T e^{-r(t)}\left(\langle\,X(s)-G(\phi(s)),\phi(s)\,\rangle +\langle\,  G(\phi(s)),\tilde{u}(s)\,\rangle\right)ds\nonumber\\&&+\int_0^T\int_Ze^{-r(s)}\left(2(Y(s),g(s,\phi(s),z))-\|g(s,
\phi(s),z)\|^2\right)\nu(dz)ds].
 \end{eqnarray}
On the other hand, applying the It\^{o} formula to the process  $e^{r(T)}\|u(T)\|^2$, taking the real part and then taking expectation, we obtain
\begin{eqnarray}\label{25}&&\mathbb{E}\left(e^{-r(T)}\|u(T)\|^2-\|u_n(0)\|^2\right)\nonumber\\&=&-\mathbb{E}\,\int_0^T e^{-r(s)}r^{\prime}(s)\|u(s)\|^2ds+2\mathbb{E}\,\textmd{Re}\int_0^T e^{-r(s)}\langle\,X(s), \tilde{u}(s)\,\rangle ds\nonumber\\&&+\mathbb{E}\,\textmd{Re}\int_0^T\int_Ze^{-r(s)}\|Y(s,z)\|^2\eta(ds,dz).
\end{eqnarray}
Therefore, in view of (\ref{24}) and (\ref{25}), we infer that
\begin{eqnarray}&&\mathbb{E}\int_0^T (e^{-r(s)}[-r^{\prime}(s)\|\tilde{u}(s)-\phi(s)\|^2+2\,\textmd{Re}\langle X(s)-G(\phi(s)),\tilde{u}(s)-\phi(s)\rangle\nonumber\\&&+\int_Z\|Y(s,z)-g(s,\phi(s),z)\|^2\nu(dz)])ds\leq 0.\nonumber
\end{eqnarray}
Taking $\phi=\tilde{u}$, we can obtain that $Y(s,z)=g(s,\tilde{u}(s),z)$ a.s. If we set $\phi=\tilde{u}-\epsilon v,\,v\in L^\infty ([0,T]\times\Omega,V),\,\epsilon>0$, it yields that
$$\mathbb{E}\int_0^T e^{-r(s)}[ -r^{\prime}(s)\epsilon^2\|v(s)\|^2+2\epsilon\,\textmd{Re}\langle X(s)-G(\tilde{u}(s)-\epsilon v(s)),v(s)\rangle ]ds\leq 0.$$
Hence
$$\mathbb{E}\int_0^T e^{-r(s)}[ -r^{\prime}(s)\epsilon\|v(s)\|^2+2\,\textmd{Re}\langle X(s)-G(\tilde{u}(s)-\epsilon v(s)),v(s)\rangle ]ds\leq 0.$$
Since $v$ is arbitrary, letting $\epsilon\rightarrow 0$, we get $X(s)=G(\tilde{u}(s))$ a.s. The proof of existence of the solution is complete.\\
\indent Step 4 : It remains to prove the uniqueness of the solution to SPDE (\ref{3}).
Suppose $\omega(t)=u_1(t)-u_2(t)$, where $u_1,\,u_2$ are the solutions of (\ref{3}) with initial conditions $u_{1}(0),\,u_{2}(0)$, respectively. We define the stopping time:$$\tau_N:=\inf\{t\geq 0,\,\|u_1(t)\|^2\geq N\}\wedge \inf\{t\geq 0,\,\|u_2(t)\|^2\geq N\}\wedge T.$$ Applying It\^{o} formula to $e^{-r(t\wedge\tau_N)}\|\omega(t\wedge\tau_N)\|^2$, taking the real part and then taking the expectation, we have
\begin{eqnarray}
&&e^{-r(t\wedge\tau_N)}\mathbb{E}\left(\|\omega(t\wedge\tau_N)\|^2-\|u_{1}(0)-u_{2}(0)\|^2\right)\nonumber\\&=&-\mathbb{E}\,\int_0^{t\wedge\tau_N} e^{-r(s)}r^{\prime}(s)\|\omega(s)\|^2ds\nonumber\\&&+2\,\mathbb{E}\,\textmd{Re}\int_0^{t\wedge\tau_N} e^{-r(s)}\langle G(u_1(s))-G(u_2(s)), \omega(s)\rangle ds\nonumber\\&&+\mathbb{E}\,\int_0^{t\wedge\tau_N}\int_Z e^{-r(s)}\|g(s,u_1(s),z)-g(s,u_2(s),z)\|^2\nu(dz)ds.\nonumber
\end{eqnarray}
It follows from (\ref{34}) that
\begin{eqnarray*}
&&e^{-r(t\wedge\tau_N)}\mathbb{E}\left(\|\omega(t\wedge\tau_N)\|^2-\|u_{1}(0)-u_{2}(0)\|^2\right)\nonumber\\&=&\mathbb{E}\int_0^{t\wedge\tau_N}e^{-r(s)}[-r^{\prime}(s)\|u_1(s)-u_2(s)\|^2+2\,\textmd{Re}\,\langle\, G(u_1(s))-G(u_2(s)),u_1-u_2\,\rangle\\&&+\int_Z\|g(s,u_1(s),z)-g(s,u_2(s),z)\|^2\nu(dz)]ds\nonumber\\
&\leq&0\nonumber,
\end{eqnarray*}
which implies the uniqueness of the solution and puts an end of the proof of Theorem 2.5.
\\\\

\noindent\textbf{Acknowledgements}
\bigskip
\vspace{2mm}

LL is supported in part by the NSF from Jiangsu
province BK20171029 and the NSF of the Jiangsu Higher Education Committee of China No. 14KJB110016. HG is supported by a China NSF Grant No. 11531006, 11771123 and PAPD of Jiangsu Higher Education Institutions.

\bigskip
\vspace{4mm}

\bibliographystyle{model1b-num-names}

\end{document}